\documentclass{lmcs}
\pdfoutput=1

\usepackage{lastpage}
\lmcsdoi{15}{3}{7}
\lmcsheading{}{\pageref{LastPage}}{}{}%
{Oct.~01,~2018}{Jul.~29,~2019}{}

\usepackage{amsmath,amssymb}
\usepackage{graphicx}
\usepackage[numbers]{natbib}
\usepackage{url}
\usepackage{epsfig}
\usepackage{amsgen}
\usepackage{amssymb}
\usepackage{amsmath}
\RequirePackage{bussproofs}
\newcommand{\modusponens}{\mathit{modusponens}}

\title{A Forgotten Theory of Proofs?}
\author{E.Engeler}
\address{ETH Zurich}
\email{erwin.engeler@math.ethz.ch}

\begin{document}
\maketitle

The Hilbert Program in G\"ottingen was winding down in the early 1930s. By then it was mostly in the hands of Paul Bernays who was writing the first volume of \emph{Grundlagen der Mathematik}. Hermann Weyl had succeeded David Hilbert. There were three outstanding doctoral students in logic: Haskell B.Curry, Saunders MacLane and Gerhard Gentzen.\footnote{As Bernays was only Dozent (and was soon to be dismissed as foreign, Swiss and of jewish ancestry), Weyl was official thesis advisor, and of course took personal interest until he also left (for Princeton, Bernays for ETH.)} These three students are at the beginning of three threads in mathematical logic: Combinatory Logic (Curry), Proof Theory (Gentzen) and Algebra of Proofs (MacLane), the last one essentially forgotten, except perhaps for some technical results useful in computer algebra (cf. Newman's Lemma).

The present author, also a student of Bernays, when looking up this mathematical ancestry, was fascinated by the contrast between MacLane's enthusiasm about the ideas in his thesis as expressed in his letters home\footnote{Excerpts in \citep{saunders:2005}}, and the almost complete absence of any mathematical follow-up. Is it possible that mathematical development has passed something by, just because MacLane did not find resonance for this work and, back in the U.S., was soon successful, as the strong mathematician he was, in other fields. His work on the conceptual structure of mathematics, category theory and its pervading influence throughout mathematics, is well known~\citep{mclarty:2007}.


In this essay we attempt to revive the idea of an algebra of proofs, place MacLane's thesis work and its vision in a new framework and try to address his original vision. He believed that the mathematical enterprise as a whole should be cast into a grand algebraic framework. If completed, it would embody axiomatic, methodology of inventing and proving, and would underly mathematical exposition; a youthful dream.

\section{Prolegomena to an algebra of mathematical thoughts}

 Let us first talk about thinking. Thinking means to apply thoughts to thoughts, thoughts being things like concepts, impressions, memories, activities, projects -- anything that you can think about, including mathematics. And, of course, the results of applying a thought to a thought. Thinking is free, all combinations of thoughts are admitted into the universe of thoughts. As a mathematician I perceive here the structure of an algebra: Thoughts are the elements of the algebra and applying a thought $X$ to a thought $Y$ is a binary operation which results in the element $X \cdot Y$, again a thought.


Mathematical thoughts are about sets of \emph{definitions, problems, theorems, proofs and proof-strategies}. In the present context, to do mathematics means to make a selection of such sets, states of knowledge and proof procedures as it were, and apply these sets to each other. To mathematize this idea, we need to represent states of mathematical knowledge and the pursuit of its development in a form that permits an application operation between them. Let us first experiment with formalized mathematics and its states of knowledge.


Mathematical logic aims to represent mathematics by a system based on a formal language. Formal mathematical thoughts thereby consist of sets of statements (axioms, theorems) and proof-trees. Take propositional logic. Let $A$ be the set of propositional formulas $a, b, c, \dots$ composed from some atomic propositions by some connectives such as $\land, \lor, \supset, \lnot$. A formal proof has the form of a tree such as

\begin{prooftree}
\AxiomC{a}   \AxiomC{b}
\BinaryInfC{c}
	\AxiomC{d}
	\BinaryInfC{g}
					\AxiomC{e}   \AxiomC{f}
					\BinaryInfC{h}
				\BinaryInfC{k}
\end{prooftree}

In an obvious notation, this tree would be rendered as

\[
\{\{\{a,b\} \rightarrow c, d \} \rightarrow g, \{e,f\} \rightarrow h \} \rightarrow k .
\]

Such a proof  can be \emph{parsed} differently in order to reflect the conceptual structure of the proof -- which in fact originally may have progressed through the development or employment of various auxiliary theorems and general lemmas. For example, $g$ may be a lemma and the proof of $k$ starts with this lemma and $e$ and $f$:

\[
\{\{\{a,b\} \rightarrow c, d \} \rightarrow g\} \rightarrow (\{\{e,f\} \rightarrow h \} \rightarrow k).
\]

Another parsing would be:

\[
\{\{\{a,b\} \rightarrow c \} \rightarrow (\{d\} \rightarrow g) \} \rightarrow (\{\{e,f\} \rightarrow h \} \rightarrow k ).
\]

All of these denote the same tree but give different narratives/ interpretations of the formal proof. They are what will be called \emph{proof-expressions} and denoted by lower-case letters such as $x$, $y$, $z$ from the latter part of the alphabet. The set $P$ of proof-expressions is built up recursively from $A$:

\begin{align*}
P_0 &= A
\\
P_{n+1} &= P_n \cup \{ \alpha \rightarrow x : x \in P_n,  \alpha \subseteq P_n\ \text{finite}\}
\\
P &= \bigcup_n P_n.
\end{align*}

Of course, these \textquotedblleft proof-expressions\textquotedblright\, represent formal proofs only in the case that the arrows correspond to legal steps in a formal proof (here of propositional logic); of this later.
The result of the proof denoted by a proof expression $x$ is $x$ itself if it is a propositional formula, an element of $A$; otherwise, if $x$ is composite it is the root of the tree denoted by $x$. We let $x^{\vdash}$ denote this root, it is a propositional formula.

Sets of proof-expressions are denoted by capital letters $X$, $Y$, \dots{} or by special symbols introduced as cases arise. Such sets represent \textquotedblleft mathematical thoughts\textquotedblright\, in the sense of the introduction to this section, here restricted to the realm of formal propositions.-- To complete the picture there, it remains to specify the operation of application, $X \cdot Y$ as follows:

\[
X \cdot Y  =  \{ x : \exists \alpha \subseteq Y, \alpha \rightarrow x \in X \}.
\]

This definition is best understood if $X$ is considered as a sort of graph of a (partial and many-valued) function, each of its elements $\alpha \rightarrow x$ associating an argument(-set) $\alpha$ to a value $x$. By this operation the set of subsets of $P$, i.e. the set of mathematical thoughts, becomes an algebraic structure, the algebra $\mathcal{P}$ of propositional thoughts.

\emph{Modus Ponens} is the thought which applied to the set of formulas $\{ a \supset b, a\}$ produces $b$. Correspondingly, $[\modusponens]$ as an element of the algebra $\mathcal{P}$ contains at least the one element $\{a \supset b, a\} \rightarrow b$; we posit that it consists of all elements of that form.
Thus, if $X$ is a set of propositional formulas, $ [\modusponens] \cdot X$ is the set of all propositional formulas provable from the set of propositional statements in $X$ in one step. Compare this with the usual notation

\[
\frac{a \supset b \qquad b} {\quad b} [\modusponens],
\]
specifying the proof-rule on the right.


More to the point, Modus Ponens can also function as a \emph{proof-constructor}. The corresponding element of $\mathcal{P}$ is

\[
[\mathrm{MP}] = \{\{x,y\} \rightarrow (\{x, y\} \rightarrow b) : \exists a \exists b\in A\ \text{such that}\ x^{\vdash} = a \supset b, \quad y^{\vdash} = b\}.
\]

$[\mathrm{MP}]\cdot X$ combines \emph{proofs} of formulas  $a \supset b$ and $a$ to a proof of $b$. An iteration starting with $X_0 = X$ using $[\mathrm{MP}]$ produces the propositional theory $T$ of $X$ restricted to the one proof-rule: $ X_{i+1} = X_i \cup [\mathrm{MP}] X_i$,

And so on, to develop propositional logic as an algebraic proof-theory of $\mathcal{P}$, a tentative partial realization of MacLane's \emph{Logikkalkul}.


Instead, we take another elementary example, finitely presented groups:
\\
Let $A$ be the set of terms $u, v, \dots$ built up from variables and a finite set of constants (\textquotedblleft generators\textquotedblright\,) denoting some elements of a group $G$ by the operations of multiplication, inverse and the unit element. Finite sets of constant terms, called relations, constitute a group-presentation. Based on $A$ we construct a \emph{calculus of reductions} $\mathcal{R}$ starting from a set of reduction-expressions $x, y, \dots$ analogously to $P$ above, (most of which of course would not denote valid reductions). Valid reductions are based on group laws such as associativity and on the relations given by the presentation:

The associative law, when applied to a reduction-expression $x$, replaces a sub-term of the final term $x^{\vdash}$, assuming it has the form $u (v w)$, by $(u v) w$, or $(u v) w$ by $u (v w)$. Let $[\mathrm{ASS}]$ denote this element of $\mathcal{R}$, hence $[\mathrm{ASS}]$ is the set of all  $ \{x\} \rightarrow t $, where $t$ results from $ x^{\vdash}$ by  substituting  some  sub-term $u(vw)$ or $(uv)w$ of $x^{\vdash}$  by  $1$. Similarly for inverse law: $[\mathrm{INV}]$  replaces sub-terms $u u^{-1}$  or $u^{-1} u$  of  $x^{\vdash}$ by  $1$. The identity law is realized as an operation $[\mathrm{ID}]$ on reductions, using replacements of $u 1$ or $1 u$ by $u$.

Relations $r_1, \dots r_n$ of the presentation give rise to reduction laws and therefore to reduction- constructors $[r_i]$. For example, if $ r_1 = g_1 g_2^{-1} g_1$ with generators $g_1, g_2$, then $[r_1]$ is the set  of all $ \{x\} \rightarrow t$, where $t$\ results from $x^{\vdash}$  by substituting some sub-term $g_1 g_2^{-1} g_1$ by  $1$ .

Example: To construct a reduction (by \textquotedblleft normalization\textquotedblright\,) of the term $(s t^{-1}) t $ we start with the set $X$, consisting of this term, and use the three operators in succession, resulting in a linear reduction-tree $x$ with $ x^{\vdash}$ equal to $1$:

\[
[\mathrm{ID}] \cdot ([\mathrm{INV}] \cdot ([\mathrm{ASS}] \cdot \{(s t^{-1}) t \}))^{\vdash} = 1.
\]

Taking the closure of $[\mathrm{ASS}] \cup [\mathrm{INV}] \cup [\mathrm{ID}] \cup [r_1] \dots \cup [r_n]$ under iteration as above, we obtain an object $[\mathrm{ALG}]$ of $\mathcal{R}$ which, applied to $X$ gives its normalization, $[\mathrm{ALG}] \cdot X$ in this finitely generated group.

\section{MacLane's Thesis and its Vision, Revisited}

The above example is from MacLane's thesis \emph{\textquotedblleft Abgek\"{u}rzte Beweise im Logikkalkul\textquotedblright\,}.\footnote{G\"{o}ttingen, Huber \& Co.1934. Reprinted in \citet{kaplansky:1979}, p.1 - 62.}  It is \textquotedblleft abgek\"{u}rzt\textquotedblright\,, shortened -- but more importantly it is a proof-template,  a formal object in a proof-manipulating system for elementary group theory, a \textquotedblleft Reduktionsbeweis\textquotedblright\,. In the original, it reads:\footnote{In the order of applications, reversed from the operational notation above; ``Th'' denotes the equation to be derived, ``Sub (4)'' denotes the application of associativity $[\mathrm{ASS}]$, etc.}

 \begin{center}
 Anfang Th, Sub (4), Sub (2), Ende (3).
 \end{center}


The Logikkalkul of MacLane takes its examples from the formal logic of \emph{Principia Mathematica}~\citep{russell-whitehead:1927}. \\
The main technical development in the thesis shows how MacLane convinces himself that his approach to proof theory suffices to treat all of mathematical logic, whose main corpus at that time was \emph{Principia}.\footnote{As it was for G\"{o}del three years before, (he strangely is not mentioned in the thesis).} \\
The basic insight is that proofs are built up from individual proof operations by composition. Correspondingly, short descriptions of logic-proofs use operators corresponding to the familiar introduction/elimination of logical connectives.

By introducing names for proofs of auxiliary theorems MacLane enriches the totality of proof-operators by names for proof-plans. It is clear that he develops the rudiments of a calculus of such expressions for proof-operators.

Admittedly, MacLane's algebraisation of deduction processes at first sight does not look very impressive. Today; but to actually complete the project, there were tedious and occasionally delicate technical details of substitution, replacement etc. to be handled. In fact, what MacLane did was at the start of a mathematics of symbol manipulations systems which later became computer algebra and computational logic, (cf. normal forms, confluence, etc.). Later in life, MacLane was aware of this~\citep{maclane:1979}.

But in 1932 MacLane's aim was broader, it was to study all formal and informal proof activities as a mathematical subject. \textquotedblleft  ... great ideas and magnificent generalizations, ...suggest all sorts of vague possibilities of applications to other fields...\textquotedblright\,\footnote{Quote in~\citep{saunders:2005}, p.59}. He submitted these ideas in a first draft of his thesis to Weyl who apparently found it too ambitious a project. Upon the wise counsel of Bernays he then confined himself to \textquotedblleft the kind of ideas [of] the professors here \textquotedblright\,\footnote{Ibid.p.60}, namely the formal manipulation and specifically the simplification of proofs. This is in fact the recently unearthed 24th problem that Hilbert had prepared for the famous 1900 list~\citep{thiele:2003}!
Some of the original targets remain as remarks in the thesis, in particular a formal analysis of the beautiful lectures of Weyl with its transparent non-formal proofs, a book by E.H.Moore, and intuitionism.


But now I'm puzzled.


\subsection{Puzzle: Relation to Curry}
As shown above, proofs, including proof plans, can be viewed as algebraic objects with an operation of composition. These form an algebraic structure which is in fact a \emph{combinatory algebra}, a model of Curry's combinatory logic. Moreover, the formation of arbitrary proofs by combination of proofs conform to the basic axiom scheme of  combinatory logic. Curry was a \textquotedblleft Kommilitone\textquotedblright\, ( roughly: a fellow-student) of MacLane, whom he remembered in his autobiography as ``a good friend of mine from G\"{o}ttingen''. He is not mentioned in the thesis.  Had MacLane presented his \textquotedblleft Logikkalkul\textquotedblright\, in the form proposed above, he would have found a model, and therefore a consistency proof, of combinatory logic just by looking at the structure $\mathcal{P}$, the elementary algebraic proof system of section 1. Whether he would have constructed the combinators $\boldsymbol{S}$  and $\boldsymbol {K}$ as explicit objects in $\mathcal{P}$ is questionable. This had to wait almost fifty years to this authors construction~\citep{engeler:1981},  which bases the Plotkin-Scott model of the Lambda Calculus on arbitrary set (useful for applications to non-numeric modeling interactive systems, logic programming,  and more recently to neuroscience.)  $\boldsymbol K$ and $\boldsymbol S$  are examples of a general algorithm that compiles expressions $\phi(X_1, \dots X_n)$ to yield an explicit "combinator" $[\mathrm{phi}]$ with $(( \dots ([\mathrm{phi}] \cdot X_1)\cdot X_2) \dots \cdot X_n) = \phi(X_1, \dots X_n)$.
\\
Applied to $\phi(X,Y) = X$, respectively $\phi(X,Y,Z) = (X Y ) (X Z)$ this compilation yields

\[
\boldsymbol{K} = \{ \{y\} \rightarrow ( \emptyset \rightarrow y) : y \in P\},
\] and
\[
\boldsymbol{S} = \{ \{ \tau \rightarrow ( \{r_1, \dots, r_n\} \rightarrow s)\} \rightarrow (\{\sigma_1 \rightarrow r_1, \dots, \sigma _n \rightarrow r_n\} \rightarrow ( \sigma \rightarrow s)):
\]
\[
n \geq 0, r_1, \dots r_n \in P, \tau \cup _i \sigma_i = \sigma \subseteq P, \sigma \quad \text{finite}\} .
\]

\subsection{Puzzle: Relation to Gentzen}

Gerhard Gentzen was more than just MacLane's contemporary at G\"{o}ttingen; indeed he translated MacLane's thesis from his English into the required German, at least in part.\footnote{Letter of MacLane to Menzler~\citep{menzler-trott:2007}.}  And he worked on his own famous thesis at just this time, in which he does also treat of normalisation and operations on proofs.  What did MacLane know or foresee and communicate on all this with Gentzen, (who is not mentioned in the thesis)?

In fact Gentzen was the one person who could have furthered MacLane's original program; if he only would have made the connection to Church's Lamba Calculus, published at the period, and cast his work in that framework. This was finally accomplished by \citet{barendregt-ghilezan:2007} receiving the predicate \textquotedblleft Theoretical Pearl\textquotedblright\,.


There may be a personal resolution of this puzzle: MacLane, in the last months of his stay in Germany was freshly married and looked toward home eagerly, and with some relief from his exposure to the frightening rise of Nazism and its disruptive incursion into academia.

The present reconsideration of his youthful dream may perhaps revive some of its parts in view of the generality of our notion of \textquotedblleft mathematical thoughts\textquotedblright\, and of their compositional structure and its algebra.

\bibliographystyle{plainnat}
\bibliography{references}

\section*{Dedication}
This is to remember my good old friend Corrado Boehm.

Our mathematical beginnings overlapped at the ETH in Zurich: Corrado was an assistant of Prof.E.Stiefel in applied mathematics in the late 1940s. As such he was sent to Bavaria to inspect the hidden Zuse Z4 computer in view of its operability. With him was Ambros Speiser, who went on to found the Zurich IBM Research Laboratory (which in his time produced four Nobel prizes; he was one of the first presidents of IFIP). The machine was then smuggled into Switzerland and became the first programmable scientific computer on the continent.

We must have met around Z4 and in Bernays' logic seminar in the early 1950's. With Richard Buechi, we were another triplet of his doctoral students.  We diverged somewhat, I went into model theory and universal algebra and the USA, while he stayed with theoretical computer science and returned to Italy. Corrado and I connected again on the semantics of programming languages and, of course, on Combinatory Logic.

\end{document}